\documentclass[12pt]{article}
\usepackage{a4}

\DeclareMathAlphabet{\eufrak}{U}{}{}{}  
\SetMathAlphabet\eufrak{normal}{U}{euf}{m}{n}
\SetMathAlphabet\eufrak{bold}{U}{euf}{b}{n}

\newcommand{\eins}{{\rm 1\!\!1}}
\newcommand{\nz}{{\rm I\!N}}
\newcommand{\gz}{{\rm Z\!\!Z}}
\newcommand{\rz}{{\rm I\!R}}
\newcommand{\cz}{{\,\rm{\sf I}\!\!\!C}}

\newcommand{\gfrak}{\eufrak{g}}

\newcommand{\expec}{{\rm I\!E}}
\newcommand{\probab}{{\rm I\!P}}
\newcommand{\bende}{\hspace{\fill}\rule{2mm}{3mm}\\ \rule{0mm}{2mm}}

\newcommand{\qbinomi}[2]{\left(\begin{array}{c} #1 \\ #2
\end{array}\right)_q}

\newtheorem{lem}{Lemma}[section]
\newtheorem{prop}[lem]{Proposition}
\newtheorem{theo}[lem]{Theorem}
\newtheorem{cor}[lem]{Corollary}
\newtheorem{defin}[lem]{Definition}

\begin{document}
\title{L\'evy Processes on $U_q(\gfrak)$ as Infinitely Divisible Representations\thanks{ASI-TPA/13/99 (TU Clausthal);\ 6/99 (Preprint-Reihe Mathmatik, 
Univ. Greifswald)}\phantom{.}\thanks{This work was prepared in a cooperation supported by Procope.}\phantom{.}\thanks{V.K.D. is thankful also to the Deutsche Forschungsgemeinschaft (DFG) for the financial  support as Guest Professor at ASI (TU Clausthal) and to the Bulgarian National Research Foundation for partial support under contract $\Phi$-643.}}
\author{V.K. Dobrev\thanks{Arnold Sommerfeld Institut f\"ur Mathematische Physik, Technische Universit\"at Clausthal, Leibnizstra{\ss}e 10, D-38678 Clausthal-Zellerfeld, Germanny, permanent address: Bulgarian Academy of Sciences, Institute of Nuclear Research and Nuclear Energy, 72 Tsarigradsko Chaussee, 1784 Sofia, Bulgaria},\ H.-D. Doebner\thanks{Arnold Sommerfeld Institut f\"ur Mathematische Physik, Technische Universit\"at Clausthal, Leibnizstra{\ss}e 10, D-38678 Clausthal-Zellerfeld, Germanny},\  U. Franz\thanks{Institut f\"ur Mathematik und Informatik, Ernst-Moritz-Arndt-Universit\"at Greifswald, Friedrich-Ludwig-Jahn-Stra{\ss}e 15a, D-17487 Greifswald, Germany}, and\ R. Schott\thanks{IECN and LORIA, Universit\'e Henri Poincar\'e-Nancy 1, B.P.\ 239, F-54506 Vand{\oe}uvre-l\`es-Nancy, France}}
\date{}
\maketitle

\begin{abstract}
L\'evy processes on bialgebras are families of infinitely divisible representations. We classify the generators of L\'evy processes on the compact forms of the quantum algebras $U_q(\gfrak)$, where $\gfrak$ is a simple Lie algebra. Then we show how the processes themselves can be reconstructed from their generators and study several classical stochastic processes that can be associated to these processes.
\end{abstract}

\section{Introduction}

L\'evy processes on involutive bialgebras made their first appearance in a model of the laser studied by von Waldenfels, cf.\ \cite{waldenfels84}. Their algebraic framework was formulated in \cite{accardi+schuermann+waldenfels88}, for their general theory see \cite{schuermann93}.

Let ${\cal B}$ be an involutive bialgebra, i.e.\ an involutive unital associative algebra over $\cz$ with two $*$-homomorphisms $\Delta:{\cal B}\to{\cal B}\otimes{\cal B}$ and $\varepsilon:{\cal B}\to \cz$ satisfying
\begin{eqnarray*}
(\Delta\otimes{\rm id})\circ\Delta &=& ({\rm id}\otimes \Delta)\circ \Delta \\
(\varepsilon\otimes {\rm id})\circ \Delta &=& ({\rm id}\otimes \varepsilon)\circ\Delta = {\rm id},
\end{eqnarray*}
and let $\pi_i=(\pi_i,H_i,\Omega_i)$ be (cyclic) representations of ${\cal B}$ on some pre-Hilbert spaces $H_i$ with the vacuum vectors $\Omega_i\in H_i$, $i=1,\ldots,n$. Then the product of these representations is the representation
\[
\Pi_{i=1}^n \pi_i=((\pi_1\otimes\cdots\otimes \pi_n)\circ \Delta^{(n)}, V_1\otimes \cdots\otimes V_n,\Omega_1\otimes\cdots\otimes \Omega_n),
\]
where $\Delta^{(n)}$ is defined by $\Delta^{(1)}={\rm id}$, $\Delta^{(2)}=\Delta$ and $\Delta^{(n)}= (\Delta\otimes {\rm id}^{\otimes (n-2)})\circ \Delta^{(n-1)}$ for $n\ge3$. A representation $\pi$ is called infinitely divisible, if for any integer $n\ge 1$ there exists a representation $\pi^{(n)}$ such that
\[
\pi\cong \Pi_{i=1}^n \pi_i, \quad\mbox{ with }\quad \pi_i=\pi^{(n)} \quad \mbox{ for } i=1,\ldots,n.
\]
Two representations $\pi=(\pi,H,\Omega)$ and $\pi'=(\pi',H',\Omega')$ are considered as equivalent, if the vacuum expectations coincide, i.e.\
\[
\pi\cong \pi' \quad \Leftrightarrow \quad \langle \Omega, \pi(a)\Omega\rangle =  \langle \Omega', \pi'(a)\Omega'\rangle \quad\mbox{ for all } a\in{\cal B}.
\]
All L\'evy processes on $*$-bialgebras (see Definition \ref{def levy}) define examples of infinitely divisible representations of these $*$-bialgebras, see \cite{schuermann90} and also below.

L\'evy processes on $*$-bialgebras have also been investigated in relation with probability theory. A class of stochastic processes with rather surprising properties, the so-called Az\'ema martingales, have been shown to be classical versions (see Section \ref{classical}) of L\'evy processes on a non-commutative, non-cocommutative involutive bialgebra, see \cite{schuermann93} and the references therein.

In this paper we will be interested in L\'evy processes on the compact forms $\cal U$ of the Drinfeld-Jimbo quantum enveloping algebras $U_q(\gfrak)$ corresponding to the simple Lie algebras $\gfrak$ in the formulation of Jimbo, cf.\ \cite{jimbo85}. In Section \ref{prelim}, we recall the definition of these $*$-bialgebras as well as the definition of L\'evy processes and some of their elementary properties. Section \ref{classification} contains the first main result of this paper, the characterisation of the L\'evy processes on $\cal U$ in terms of their generators. It turns out that all generators are of the form
\[
\psi(u) = \langle \Omega, (\rho(u)-\varepsilon(u))\Omega \rangle, \qquad
\mbox{ for all } u\in{\cal U},
\]
where $\rho$ is a unitary representation of $\cal U$ on some pre-Hilbert
space $D$, $\Omega$ some vector in $D$ and $\varepsilon$ the counit of $\cal U$.

In Section \ref{construction}, we show how to reconstruct a realization of the process on a Bose-Fock space using quantum stochastic calculus. With the help of the explicit expressions for the Cartan elements given in Proposition \ref{group-like}, we can give a classical stochastic process in Section \ref{classical}, whose joint moments coincide with the vacuum expectation of the restriction of the process to the Cartan subalgebra (see Theorem \ref{cartan}). This process is seen to be a Poisson jump process on the lattice generated by the weights of the unitary irreducible representations of $\cal U$. We also give several other elements of $\cal U$, whose vacuum expectations can be characterised by classical processes.

\section{Preliminaries}\label{prelim}

\subsection{The Hopf algebras $U_q(\gfrak)$}

Let $\gfrak$ be any complex simple Lie algebra and 
$(a_{ij})_{1\le i,j\le n}$ be its Cartan matrix. Let $(d_i)_{1\le i
\le n}$ be non-zero integers such that $d_i a_{ij} = d_j a_{ji}$ and
the greatest common divisor of the $d_i$'s is $1$. Let furthermore 
$q\not=0$ be a complex number such that $q^{2d_i}\not= 1$ for all $i$. 
The quantum enveloping algebra $U_q(\gfrak)$  is defined \cite{jimbo85}
as the Hopf algebra generated by $e_i$, $f_i$, $k_i$, and $k^{-1}_i$, 
$i=1,\ldots,n$ with the relations
\begin{eqnarray*}
k_ik^{-1}_i = k_i^{-1} k_i =1, && k_ik_j=k_jk_i, \\
k_i e_j = q_i^{a_{ij}/2} e_j k_i, && k_i f_j = q_i^{-a_{ij}/2} f_j
k_i,\\
e_i f_j - f_j e_i &=& \delta_{ij} \frac{k_i^2 -k_i^{-2}}{q_i-q_i^{-1}},
\\
\sum_{n=0}^{1-a_{ij}}\left(\begin{array}{c}{1-a_{ij}} \\ {n}
\end{array}\right)_{q_i} (-1)^n e_i^{1-a_{ij}-n} e_j e_i^n =0 && 
\mbox{ for } i\not=j \\
\sum_{n=0}^{1-a_{ij}} \left(\begin{array}{c}{1-a_{ij}} \\ {n}
\end{array}\right)_{q_i} (-1)^n f_i^{1-a_{ij}-n} f_j f_i^n =0 && 
\mbox{ for } i\not=j \\
\Delta(e_i) = e_i\otimes k_i^{-1} + k_i\otimes e_i,\quad \Delta(f_i) &=& 
f_i\otimes k_i^{-1} + k_i\otimes g_i,\quad \Delta(k_i^{\pm 1}) = 
k_i^{\pm 1}\otimes k_i^{\pm 1}, \\
\varepsilon(e_i)=\varepsilon(f_i) = 0, && \varepsilon(k_i^{\pm 1}) =
1,\\
S(e_i) = -q_i^{-1} e_i, && S(f_i) =- q_i f_i, \quad S(k_i) = k_i^{-1},
\end{eqnarray*}
where $q_i=q^{d_i}$ and $\qbinomi{n}{\nu}$ is defined by
\[
\qbinomi{n}{\nu} = \frac{[n]_q!}{[\nu]_q! [n-\nu]_q!}, \quad [n]_q! = 
[1]_q [2]_q \cdots [n]_q, \quad [n]_q=\frac{q^n-q^{-n}}{q-q^{-1}}.
\]

We restrict ourselves to $q\in\rz\backslash\{0\}$ and also define an anti-involution
on $U_q(\gfrak)$ which is given  on the generators by
\[
(e_i)^* = f_i, \quad (f_i)^* = e_i, \quad (k_i)^*=k_i,
\]
The implementation of this anti-involution produces the quantum enveloping 
algebras $U_q(\gfrak_c)$ corresponding to the simple compact Lie algebras $\gfrak_c\,$.  
We shall use the notation: ~${\cal U}\equiv U_q(\gfrak_c)$.

\subsection{L\'evy processes on bialgebras}

We recall the definition of L\'evy processes on $*$-bialgebras, cf.\ \cite{schuermann93}.

\begin{defin}\label{def levy}
A family of $*$-homomorphisms $(j_{st})_{0\le s\le t}$ defined on a
$*$-bialgebra ${\cal B}$ with values in another $*$-algebra ${\cal A}$
with some fixed state $\Phi:{\cal A}\to \cz$ is called a L\'evy process
(w.r.t.\ $\Phi$), if the following conditions are satisfied:
\begin{description}
\item[(i)]
the images corresponding to disjoint time intervals commute, i.e.\
$[j_{st}({\cal B}),j_{s't'}({\cal B})]=\{0\}$ for $0\le s\le t\le s'\le
t'$, and expectations corresponding to disjoint time intervals
factorize, i.e.\
\[
\Phi(j_{s_1t_1}(b_1) \cdots j_{s_nt_n}(b_n)) =
\Phi(j_{s_1t_1}(b_1))\cdots \Phi(j_{s_nt_n}(b_n)),
\]
for all $n\in\nz$, $b_1,\ldots,b_n\in{\cal B}$ and $0\le s_1\le
t_1\le\cdots\le t_n$;
\item[(ii)]
$m_{\cal A}\circ(j_{st}\otimes j_{tu})\circ \Delta = j_{su}$ for all $0\le s\le t\le
u$;
\item[(iii)]
the functionals $\varphi_{st}=\Phi\circ j_{st}:{\cal B}\to \cz$ depend
only on $t-s$;
\item[(iv)] $\lim_{t\searrow s} j_{st}(b) = j_{ss}(b) = \varepsilon(b)
1_{\cal A}$ for all $b\in {\cal B}$.
\end{description}
\end{defin}

For a detailed exposition of the theory of these processes see \cite{schuermann93}, for a more accessible first introduction see also \cite[Chapter VII]{meyer95} or \cite{schuermann91}.

The functionals $\varphi_{t-s} = \varphi_{st}$ then form a convolution
semi-group of states and there exists a hermitian conditionally positive
(i.e.\ positive on ${\rm ker}\,\varepsilon$) linear functional $\psi$
such that $\varphi_t=\exp_\star t\psi = \varepsilon +t\psi +
\frac{t^2}{2} \psi\star\psi + \cdots + \frac{t^n}{n!} \psi^{\star n} +
\cdots$ and $\psi(1)=0$. Conversely, for every hermitian conditionally
positive linear functional $\psi:{\cal B}\to \cz$ with $\psi(1)=0$ there
exists a unique convolution semi-group of states $(\varphi_t)_{t\in\rz_+}$ and a unique (up to
equivalence) L\'evy process $(j_{st})_{0\le s\le t}$. The functional $\psi$ is called the generator of the L\'evy process $(j_{st})_{0\le s\le t}$.

\begin{defin}
Two quantum stochastic processes (i.e.\ families of $*$-homomorphisms) $(j^{(1)}_{st}:{\cal B}\to{\cal A}_1)$ and 
$(j^{(2)}_{st}:{\cal B}\to{\cal A}_2)$ on the same $*$-bialgebra ${\cal
B}$ are called equivalent (with respect to two fixed states $\Phi_1$ and
$\Phi_2$ on ${\cal A}_1$ and ${\cal A}_2$, respectively) if their joint
moments agree, i.e.\ if
\[
\Phi_1\left(j^{(1)}_{s_1t_1}(b_1)\cdots j^{(1)}_{s_nt_n}(b_n)\right) =
\Phi_2\left(j^{(2)}_{s_1t_1}(b_1)\cdots j^{(2)}_{s_nt_n}(b_n)\right),
\]
for all $n\in\nz$, $b_1,\ldots,b_n\in{\cal B}$ and $s_1,\ldots,
s_n,t_1,\ldots,t_n\in\rz_+$.
\end{defin}

Let us now show that L\'evy processes give indeed infinitely divisible representations, as we claimed in the Introduction. Let $\pi_{st}=(\pi_{st},H_{st},\Omega_{st})$ be the GNS representation induced by $j_{st}:{\cal B}\to{\cal A}$, i.e.\ $H_{st}$ is the pre-Hilbert space obtained by taking the quotient of $\cal B$ with respect to the null space ${\cal N}_{st}=\{a\in{\cal B}; \Phi(j_{st}(a^*a))=0\}$ with the inner product induced from the inner product $\langle a, b\rangle=\Phi(j_{st}(a^*b))$ on $\cal B$, $\pi_{st}$ is the representation of $\cal B$ on $H_{st}={\cal B}/{\cal N}_{st}$ induced from left multiplication, and $\Omega_{st}$ is the image of the unit element $1$ under the canonical projection from $\cal B$ to $H_{st}$. Then, by property (iii) of Definition \ref{def levy}, $(\pi_{st}:{\cal B}\to L(H_{st}))$ is equivalent to $(j_{s't'}:{\cal B}\to{\cal A})$ (with respect to the states $\langle\Omega_{st},\pi_{st}(\cdot)\Omega_{st}\rangle$ and $\Phi$ on $L(H_{st})$ and!
 $\cal A$, respectively) if the 
intervals $(s,t)$ and $(s't')$ have the same lenght. The product of such representations $\pi_{s_1t_1},\ldots,\pi_{s_nt_n}$ is therefore by property (ii) equivalent to $j_{s_1+\cdots+s_n,t_1+\cdots+t_n}\cong \pi_{s_1+\cdots+s_n,t_1+\cdots+t_n}$. Therefore, if we want to write $j_{st}\cong\pi_{st}$ as an $n$-fold product
\[
\pi_{st}\cong \Pi_{i=1}^n \pi_i, \quad\mbox{ with }\quad \pi_i=\pi^{(n)} \quad \mbox{ for } i=1,\ldots,n,
\]
it is sufficient to take $\pi^{(n)}=\pi_{(t-s)/n}$. This proves that $\pi_{st}$ is indeed infinitely divisible.

For a given generator $\psi$ one can construct the so-called Sch\"urmann
triple $(\rho,\eta,\psi)$ consisting of a unitary representation $\rho$
of $\cal B$ on some pre-Hilbert space $D$, a $(\rho,\eta)$-1-cocycle
$\eta:{\cal B}\to D$, and the generator $\psi:{\cal B}\to C$ itself,
such that the following relations
\begin{eqnarray*}
\eta(ab) &=& \rho(a)\eta(b) + \eta(a) \varepsilon(b), \\
\langle a,b\rangle &=& - \varepsilon(a^*)\psi(b) +\psi(a^*b) - \psi(a^*)
\varepsilon(b),
\end{eqnarray*}
hold for all $a,b\in {\cal B}$.

This construction goes as follows. First we define a sesqui-linear form
on ${\cal B}$ by
\[
\langle a,b\rangle_\psi =
\psi((a-\varepsilon(a)1)^*(b-\varepsilon(b)1)),
\]
for $a,b\in {\cal B}$. Since $a\to a-\varepsilon(a)1$ is a projection
from $\cal B$ to ${\rm ker}\,\varepsilon$ and since $\psi$ is positive
on ${\rm ker}\,\varepsilon$, this form is positive semi-definite. If we
quotient $\cal B$ by the nullspace of this form,
\[
{\cal N}_\psi = \{ b\in {\cal B}; \langle a,a\rangle_\psi=0\},
\]
then we obtain a pre-Hilbert space $D={\cal B}/{\cal N}_\psi$, this will
be the space on which the representation $\rho$ acts. The cocycle $\eta$
is just the canonical projection from $\cal B$ onto $D$, and the inner
product on $D$ and the sesqui-linear form on ${\cal B}$ are related via
\[
\langle \eta(a),\eta(b)\rangle = \langle a, b\rangle_\psi =
-\varepsilon(a^*)\psi(b) + \psi(a^*b) - \psi(a^*) \varepsilon(b),
\]
for $a,b\in {\cal B}$.
Since ${\cal N}_\psi\cap {\rm ker}\, \varepsilon$ is invariant under
left multiplication of elements of ${\rm ker}\,\varepsilon$, we have an
action of ${\rm ker}\,\varepsilon$ on ${\rm ker}\,\varepsilon/({\cal
N}_\psi\cap{\rm ker}\,\varepsilon) = D$ such that $\rho(a) \eta(b) =
\eta(ab)$ for all $a,b\in{\rm ker}\,\varepsilon$. This representation
can be extended to a representation of $\cal B$ on $D$, if we set
\[
\rho(a)\eta(b) = \eta(ab) - \eta(a)\varepsilon(b),
\]
for $a,b\in{\cal B}$. If we require the cocycle $\eta$ to be onto, then
the Sch{\"u}rmann triple is unique up to isometry.

We will do this construction in the reversed sense in Section
\ref{classification} in order to classify the generators of L\'evy processes on the $*$-Hopf
algebras ${\cal U}$ starting from the classification of their unitary
representations.

\section{Classification of the generators}\label{classification}

We will now give a complete classification of the generators of L\'evy processes on ${\cal U}$. We begin by proving a series of lemmas which we shall use to formulate our main result.

\begin{lem} \label{k2=k1}
Define the following ideals: $K_1={\rm ker}\, \varepsilon$ and
$K_2=K_1\cdot K_1={\rm span}\, \{ uv; u,v\in K_1\}$. Then $K_2=K_1$.
\end{lem}
{\bf Proof:}
 $K_1$ is generated by $e_i$, $f_i$, $k_i-1$ and $k^{-1}_i-1$. The
relation $k_ie_j = q_i^{a_{ij}/2} e_j k_i$ implies $e_j\in K_2$ since it
can be written as
\[
e_j = \frac{(k_i-1)e_j(k_i^{-1}-1) + (k_i-1)e_j +
e_j(k_i^{-1}-1)}{q_i^{a_{ij}/2}-1}
\]
if we choose $i$ such that $a_{ij}\not=0$. Similarly, we get $f_j\in
K_2$. 
And $k_1-1,k_i^{-1}-1\in K_2$ follows from
\[
(q_i-q_i^{-1}) [e_i,f_i] - (k_i-1)^2 +(k_i^{-1}-1)^2 = 2(k_i-k_i^{-1}) \in
K_2,
\]
since $k_i-1= \frac{1}{2}\left(k_i(k_i-k_i^{-1}) - (k_i-1)^2\right)$ and
$k_i^{-1}-1= -\frac{1}{2}\left(k_i^{-1}(k_i-k_i^{-1})
+(k_i^{-1}-1)^2\right)$.
\bende

\begin{lem}\label{coboundary}
Suppose we have a second order Casimir element $C$ in $\cal U$ such that
for all unitary irreducible representations $\pi$ of $\cal U$ except the
one-dimensional $\pi(C)-\varepsilon(C)$ is invertible. Let $\rho$ be an
arbitrary unitary representation of ${\cal U}$ on some pre-Hilbert space
$D$. Then all $(\rho,\varepsilon)$-1-cocycles are trivial, i.e.\ there
exists a vector $\Omega_\eta\in D$ such that
\[
\eta(u) = (\rho(u)-\varepsilon(u)) \Omega_\eta \qquad \mbox{ for all }
u\in {\cal U}.
\]
\end{lem}
{\bf Proof:}
The cocycle equation $\eta(uv)=\rho(u)\eta(v)+\eta(u)\varepsilon(v)$
implies $\eta(1)=0$, so that it is sufficient to determine $\eta$ on
$K_1={\rm ker}\,\varepsilon$.

The representation $\rho$ is direct sum of the unitary irreducible ones,
$\rho=\sum_{\lambda\in\Lambda_\rho} \pi_\lambda$,
$D=\sum_{\lambda\in\Lambda_\rho} D_\lambda$. Let $\Omega_\eta$ be the
unique vector in $D$ that satisfies
\[
\eta(C) = (\rho(C)-\varepsilon(C))\Omega_\eta,
\]
and that has no component in the one-dimensional representations. We now
have to show that $\eta$ is equal to the cocyle defined by
$\Omega_\eta$.

$uC=Cu$ for all $u\in{\cal U}$ implies
\[
\rho(u)\eta(C) + \eta(u)\varepsilon(C)= \eta(uC) = \eta(Cu) =
\rho(C)\eta(u) + \eta(C) \varepsilon(u)
\]
and thus $(\rho(C)-\varepsilon(C))\eta(u)=\rho(u)\eta(C)$ for $u\in
K_1$. Therefore $\eta(u)-\rho(u)\Omega_\eta$ has to be contained in the
one-dimensional components of $D=\sum_{\lambda\in\Lambda_\rho}
D_\lambda$ and the cocycle $\bar\eta$ defined by
$\bar\eta(u)=\eta(u)-\rho(u)\Omega_\eta$ for $u\in{\cal U}$ is an
$(\varepsilon,\varepsilon)$-1-cocycle. This implies that $\bar\eta$ is
zero on $K_2$, and, by Lemma \ref{k2=k1}, on $K_1$, and thus on all of
$\cal U$.
\bende

\begin{lem}
The generator $\psi:{\cal U} \to \cz$ is uniquely determined by $\rho$
and $\eta$.
\end{lem}
{\bf Proof:} This follows immediately from Lemma \ref{k2=k1}, since
$\psi(1)=0$ and ${\cal U}=\cz 1 \oplus K_1$.
\bende

\begin{lem}\label{casimir}
There exists a second order Casimir element $C$ in
$\cal U$ such that for all unitary irreducible representations $\pi$
of ${\cal U}$ except the one-dimensional $\pi(C)-\varepsilon(C)$ is
invertible. \end{lem}
{\bf Proof:} First we recall that there exists a  second order 
Casimir element $C_2$ in $U_q(\gfrak)$, cf. e.g., 
\cite{sklyanin83,jimbo85,jimbo86}. For all finite-dimensional highest
weight representations of  $U_q(\gfrak)$, except for the one-dimensional 
$\pi_{\rm id}\,$, we have $\pi(C_2)\neq\varepsilon(C_2)$. 
This is so since  $\pi(C_2)\neq \pi_{\rm id}(C_2) (=\varepsilon(C_2))$ 
if $\pi \neq \pi_{\rm id}$.  
The last fact follows from the results of Rosso \cite{rosso88}, 
but can also be  
obtained by supposing the inverse and then getting a contradiction 
in the limit $q\to 1$.  The same facts hold for $C_2^*$ - the image 
of $C_2$ under the involution producing ${\cal U}$, and also 
for $C=(C_2+C_2^*)/2$. Now it only remains to note that $C$ is 
invariant under the involution and thus is the required Casimir 
of ${\cal U}$, and that the unitary irreducible representations of ${\cal U}$ are in 1-to-1
correspondence with the finite-dimensional highest weight
representations of  $U_q(\gfrak)$ and are obtained from the latter 
by use of the involution.\bende

We summarize these results in the following theorem.
\begin{theo}\label{generator}
Every generator on $\cal U$ is of the form
\[
\psi(u) = \langle \Omega, (\rho(u)-\varepsilon(u))\Omega \rangle, \qquad
\mbox{ for all } u\in{\cal U},
\]
where $\rho$ is a unitary representation of $\cal U$ on some pre-Hilbert
space $D$ and $\Omega$ some vector in $D$.
\end{theo}

{\bf Remark}: The correspondence between generators $\psi$ and
Sch{\"u}rmann triples $(\rho,\eta,\psi)$ is 1-to-1, if we require $\eta$
to be surjective. Furthermore, for a given cocycle $\eta$ the choice of
$\Omega_\eta$ is unique, if we demand that $\Omega_\eta$ has no
components in the trivial one-dimensional representations. Therefore the
correspondence between L\'evy processes on $\cal U$ and triples
$(\rho,D,\Omega)$ consisting of a unitary representation $\rho$ on a
pre-Hilbert space $D$ and a vector $\Omega\in D$ becomes a bijection, if
we impose the following two conditions:
\begin{enumerate}
\item
The trivial one-dimensional representation does not appear in the direct
sum decomposition $\rho=\sum_{\lambda\in\Lambda_\rho} \pi_\lambda$,
$D=\sum_{\lambda\in\Lambda_\rho} D_\lambda$ of $(\rho,D)$.
\item
The vector $\Omega\in D$ is cyclic for $(\rho,D)$, i.e.\ $D=\rho({\cal
U})\Omega$.
\end{enumerate}

\section{Construction of the processes}\label{construction}

{}From Theorem \ref{generator} we know that every generator on $\cal U$
is given by a unitary representation $\rho$ on a pre-Hilbert space $D$
and a vector $\Omega\in D$. From the general theory follows that the
corresponding process can be constructed on the Fock space
$\Gamma(L^2(\rz_+,H))$ (or rather on a dense stable subspace thereof),
where $H$ is the Hilbert space closure of $D$, via a quantum stochastic
differential equation.

More precisely, we have the following theorem.
\begin{theo}\label{representation thm}
Let $\rho$ be a unitary representation of ${\cal U}$ on $D$ and let
$\Omega\in D$. Set $\eta(u)=(\rho(u)-\varepsilon(u))\Omega$,
$\tilde\eta(u)=\eta(u^*)$, $\psi(u)= \langle \Omega,
(\rho(u)-\varepsilon(u))\Omega \rangle$ for $u\in{\cal U}$. Then the
quantum stochastic differential equations
\[
{\rm d}j_{st} (u) = m\circ \Big(j_{st}\otimes( {\rm
d}\Lambda(\rho-\varepsilon) + {\rm d}A^*(\eta) + {\rm d}A(\tilde\eta) +
\psi {\rm d}t)\Big)(\Delta u),
\]
with the initial conditions $j_{ss}(u)=\varepsilon(u){\rm id}$ for
$u\in{\cal U}$ have solutions on a domain ${\cal
E}_D\subseteq\Gamma(L^2(\rz_+,H))$ that contains the Fock vacuum, that
is dense in the Fock space $\Gamma(L^2(\rz_+,H))$ and that is invariant
under $j_{st}({\cal U})$ for all $s\le t$, $s,t\in\rz_+$.

Furthermore, in the vacuum state $(j_{st})$ is a L\'evy process on $\cal
U$ with generator $\psi$.
\end{theo}
{\bf Proof:} The triple $(\rho,\eta,\psi)$ satisfies the conditions of
\cite[Theorem 2.3.5]{schuermann93}, therefore our theorem follows from
Sch{\"u}rmann's representation theorem (see \cite[Theorem
2.3.5]{schuermann93} and \cite[Theorem 2.5.3]{schuermann93}). For the
exact definition of the domain ${\cal E}_D$, see Page 44 of
\cite{schuermann93}.
\bende

For group-like elements (e.g.\ the $k_i$'s) we can use 
\cite[Proposition 4.1.2]{schuermann93} to get explicit expressions
without having to solve any quantum stochastic differential equations.
These expressions become particularly simple, if we act on the
exponential or coherent vectors
\[
{\cal E}(f)= \sum_{n\in\nz} \frac{f^{\otimes n}}{\sqrt{n!}}
\]
for $f\in L^2(\rz_+)\otimes D$.

\begin{prop}\label{group-like}
Let $\lambda=(\lambda_1,\ldots,\lambda_n)\in\gz^n$,
$k^\lambda=k_1^{\lambda_1}\cdots k_n^{\lambda_n}$, $f\in
L^2(\rz_+)\otimes D$, and $(j_{st})$ be the process defined in Theorem
\ref{representation thm} for the triple $(\rho,\eta,\psi)$. Then we have
\begin{eqnarray*}
&&j_{st}\left(k^\lambda\right){\cal E}(f) \\
&=& \exp\left((t-s)\psi(k^\lambda)
+ \int_s^t \langle \eta(k^\lambda), f(r)\rangle {\rm d}r\right) {\cal
E}\left( f\eins_{[0,s[\cup[t,\infty[} + (\rho(k^\lambda)f +
\eta(k^\lambda))\eins_{[s,t[}\right).
\end{eqnarray*}
In particular, in the Fock vacuum we get
\[
\varphi_{t-s}(k^\lambda) = \langle {\cal E}(0), j_{st}(k^\lambda){\cal
E}(0) \rangle = \exp\left((t-s)\psi(k^\lambda)\right).
\]
\end{prop}

\section{Classical processes}\label{classical}

In this section we will assume $q>0$.

The Cartan elements $k_i,k_i^{-1}$, $i=1,\ldots,n$, generate a
commutative sub-Hopf $*$-algebra ${\cal K}$ of $\cal U$. Therefore the
restriction of $(j_{st})$ to ${\cal K}$ is still a L\'evy process.
Furthermore, due to the self-adjointness of the $k_i$ and the
commutativity of ${\cal K}$, there exists a classical version of
$\Big(j_{st}(k_1),\ldots, j_{st}(k_n)\Big)$, i.e.\ a real-valued
stochastic process $\Big(\hat{k}_1(s,t),\ldots,\hat{k}_n(s,t)\Big)$ such
that
\begin{equation}\label{class version}
\Phi\left(j_{s_1t_1}(k^{\mu_1})\cdots j_{s_mt_m}(k^{\mu_m})\right) =
\expec\left( \hat{k}^{\mu_1}(s_1,t_1)\cdots
\hat{k}^{\mu_m}(s_m,t_m)\right)
\end{equation}
holds for all $m\in\nz$, $\mu_1,\ldots,\mu_m\in\gz^n$,
$s_1,\ldots,s_m,t_1,\ldots,t_m\in\rz_+$. Using Proposition
\ref{group-like} we can explicitly characterize the process
$\Big(\hat{k}_1(s,t),\ldots,\hat{k}_n(s,t)\Big)$.

We know that $\psi$ is of the form $\psi(u) = \langle \Omega,
(\rho(u)-\varepsilon(u))\Omega \rangle$ for all $u\in{\cal U}$ with some
unitary representation $\rho$ acting on a pre-Hilbert space $D$ and some
vector $\Omega\in D$. $D$ is a direct sum of the unitary irreducible
representation of $\cal U$ and can be decomposed into a direct sum of
eigenspaces $D=\bigoplus_\kappa E_\kappa$ of the Cartan elements
$k_1,\ldots, k_n$. Furthermore, we know that the eigenvalues are of the form
$\kappa=(q^{\lambda_1/2},\ldots,q^{\lambda_n/2})$ with
$\lambda=(\lambda_1,\ldots,\lambda_n)\in\gz^n$. Develop $\Omega$ into a
sum of eigenvectors
\[
\Omega= \sum_{\lambda\in\Lambda_\Omega} v_\lambda, \quad \mbox{ with
}\quad \rho(k_i)v_\lambda = q^{\lambda_i/2} v_\lambda
\]
and set $c_\lambda=||v_\lambda||^2$. Then we get
\[
\psi(k^\mu) = \sum_{\lambda\in\Lambda_\Omega} c_\lambda \left(q^{\langle
\mu,\lambda\rangle/2} - 1\right)
\]
for all $\mu=(\mu_1,\ldots,\mu_n)\in\gz^n$. From Proposition
\ref{group-like} we can now deduce the moments of
$\Big(\hat{k}_1(s,t),\ldots,\hat{k}_n(s,t)\Big)$, we get
\[
\Phi\Big(j_{st}(k^\mu)\Big) =
\exp\left((t-s)\sum_{\lambda\in\Lambda_\Omega} c_\lambda
\left(q^{\langle \mu,\lambda\rangle/2} - 1\right)\right)
\]
for $(\lambda_1,\ldots,\lambda_n)\in\gz^n$. We see that we can give
$\Big(\hat{k}_1(s,t),\ldots,\hat{k}_n(s,t)\Big)$ as a function of a
Poisson process on the lattice generated by the elements of
$\Lambda_\Omega$.

\begin{theo}\label{cartan}
Let $\{(N^{(\lambda)}_t)_{t\in\rz_+};\lambda\in\Lambda_\Omega\}$, be a
family of independent Poisson processes, and define a jump process
$(N_t)_{t\in\rz_+}=\Big((N_1(t),\ldots,N_n(t))\Big)_{t\in\rz_+}$ with
values in the lattice generated by the set $\Lambda_\Omega$ by
\[
N_t = \sum_{\lambda\in\Lambda_\Omega} \lambda N_{c_\lambda t},\qquad
t\in \rz_+.
\]
Then
\[
\Big(\hat{k}_1(s,t),\ldots,\hat{k}_n(s,t)\Big) = \left(
q^{(N_1(t)-N_1(s))/2}, \ldots, q^{(N_n(t)-N_n(s))/2} \right)
\]
is a classical version of $\Big(j_{st}(k_1),\ldots, j_{st}(k_n)\Big)$.
\end{theo}
{\bf Remark:} We can also give the following equivalent construction of
$(N_t)$. Set $c=\sum_{\lambda\in\Lambda_\Omega}c_\lambda$, i.e.\
$c=||\Omega||^2$. Let $(T_i)_{i\in\nz}$ be a sequence of independent,
identically distributed (i.i.d.) random variables with exponential
distribution with parameter $c$, i.e.\ $\probab(T_i< \tau) =
1-e^{-c\tau}$ for $\tau\ge 0$, and let $(\Delta_i)_{i\in\nz}$ be independent, identically distributed random variables, independent of $(T_i)$, with values in
$\Lambda_\Omega$, such that $\probab(\Delta_i = \lambda)=c_\lambda/c$.
Then we can define $(N_t)$ as
\[
N_t= \sum_{i=1}^\infty \Delta_i \eins_{\left\{\sum_{l=1}^i T_l \le
t\right\}}, \qquad t\in\rz_+.
\]
{\bf Proof:}
We have to show that Equation (\ref{class version}) is satisfied for all
$m\in\nz$, $\mu_1,\ldots,\mu_m\in\gz^n$,
$s_1,\ldots,s_m,t_1,\ldots,t_m\in\rz_+$. Without loss of generality we
can assume that $s_1\le t_1\le s_2 \le \cdots \le t_m$. Using  the
independence of the increments, we get
\begin{eqnarray*}
\Phi\Big(j_{s_1t_1}(k^{\mu_1})\cdots j_{s_mt_m}(k^{\mu_m})\Big) &=&
\Phi\Big(j_{s_1t_1}(k^{\mu_1})\Big)\cdots
\Phi\Big(j_{s_mt_m}(k^{\mu_m})\Big) \\
&=&\prod_{l=1}^m \exp\left((t_l-s_l)\sum_{\lambda\in\Lambda_\Omega}
c_\lambda \left(q^{\langle\mu_l,\lambda\rangle/2} - 1\right)\right) \\
&=& \prod_{l=1}^m \expec\left(
q^{\frac{1}{2}\sum_{\lambda\in\Lambda_\Omega}
\langle\mu_l,\lambda\rangle \left(N^{(\lambda)}_{c_\lambda t_l} -
N^{(\lambda)}_{c_\lambda s_l}\right)}\right) \\
&=& \expec\left( \hat{k}^{\lambda_1}(s_1,t_1)\cdots
\hat{k}^{\lambda_m}(s_m,t_m)\right).
\end{eqnarray*}
\bende

Comparing the right-hand-side and the left-hand-side of Equation (\ref{class version}), we get the following result.

\begin{cor}
Let ${\cal E}(0)=\sum_{\lambda\in\Lambda_{j_{st}}} v_\lambda(s,t)$ be a
decomposition of the Fock vacuum into a sum of joint eigenvectors of
$j_{st}(k_1),\ldots, j_{st}(k_n)$, i.e.\
$j_{st}(k_i)v_\lambda(s,t)=q^{\lambda_i/2} v_\lambda(s,t)$. Then the
norms of the $v_\lambda(s,t)$ are given by
\[
||v_\lambda(s,t)||^2 = \probab(N_t=\lambda), \qquad
\lambda\in\Lambda_{j_{st}}.
\]
\end{cor}

The Casimir elements also give commuting families of operators and therefore classical processes.

\begin{prop}
Let $C$ be the self-adjoint second order Casimir element.(cf.\ Lemma
\ref{casimir}). Then $(j_{0t}(C))$ has a classical version, i.e.\ there exists a real-valued stochastic process $(\hat{C}(t))_{t\in \rz_+}$ such that all joint moments of $(j_{0t}(C))$ agree with those of $(\hat{C}(t))$.

More generally, let $C_1,\ldots, C_r$ be Casimir operators, i.e.\ elements of the center ${\cal Z}({\cal U})$ of ${\cal U}$, such that $C_i^*=C_i$ for $i=1,\ldots,r$. Then there exists a classical version of $(j_{0t}(C_1),\ldots,j_{0t}(C_r))_{t\in\rz_+}$, i.e.\ an $\rz^r$-valued stochastic process with the same joint moments.
\end{prop}
{\bf Proof:} This follows from \cite[Proposition 4.2.3]{schuermann93} (see also \cite[Theorem 2.3]{franz99} for the multi-dimensional version). The commutation relations $[C_i\otimes 1,\Delta(C_j)]=0$ for $i,j=1,\ldots,r$ imply that $[j_{0s}(C_i),j_{0t}(C_j)]=j_{0s}\otimes j_{st}([C_i\otimes 1,\Delta(C_j)])=0$ for $0\le s\le t$, and therefore that the operators of the family $(j_{0t}(C_1),\ldots,j_{0t}(C_r))_{t\in\rz_+}$ commute. Since they are symmetric, their joint moments are positive and there exists a (not necessarily unique) solution of the associated moment problem. Any classical process whose distribution is given by such a solution of the moment problem is a classical version of $(j_{0t}(C_1),\ldots,j_{0t}(C_r))_{t\in\rz_+}$.
\bende

For the Lie algebra $u(n)$ and the Casimir operators
\[
G^n_m=\sum_{1\le i_1<\cdots<i_m\le n} \sum_{\pi,\sigma\in {\cal S}_m} {\rm sgn}\,(\pi\sigma) E_{i_{\pi(1)}i_{\sigma(1)}} \cdots E_{i_{\pi(m)}i_{\sigma(m)}}
\]
 of ``determinant type'', the process $({\cal G}^n_m(t)=j_{0t}(G^n_m))_{1\le m\le n, t\ge 0}$ (for a certain L\'evy process $(j_{st})$ on $U(u(n))$) has been considered by Hudson and Parthasarathy in \cite{hudson+parthasarathy94}. But a characterization of the classical process associated to this operator process is not known even for this special case, as far as we know. It would be interesting to know more about these processes. For a special process on $U_q(su(2))$ with a relatively simple generator the generator of the classical version was determined in \cite{franz99}.

We can find another element that gives us a classical process.

\begin{prop}
Let $i\in\{1,\ldots,n\}$ and set $Z=k_i^{-1}e_i+f_ik_i^{-1}$. Then
$(j_{0t}(z))$ has a classical version, i.e.\ there exists a real-valued
stochastic process $(\hat{Z}_t)$ such that
\[
\Phi\Big(j_{0t_1}(Z^{\mu_1})\cdots j_{0t_m}(Z^{\mu_m})\Big)=
\expec\left(\hat{Z}^{\mu_1}_{t_1} \cdots \hat{Z}^{\mu_m}_{t_m}\right)
\]
for all $m\in\nz$, $\mu_1,\ldots,\mu_m\in\nz$, $t_1,\ldots,t_m\in\rz_+$.
\end{prop}
{\bf Proof:} As in the preceding proposition the existence of the classical version follows from Sch\"urmann's criterium \cite[Proposition 4.2.3]{schuermann93}, since $Z\otimes 1$ and $\Delta Z=Z\otimes  k_i^{-2} + 1 \otimes Z$ commute.
\bende

\section{Outlook}

In this paper we have discussed L{\'e}vy processes on the compact 
forms $\cal U$. Further, we plan to investigate 
what is the most general class of Hopf 
algebras, on which the arguments work and analogous results can be obtained.
Naturally, we shall look first at the quantum enveloping 
algebras $U_q(\gfrak_n)$ corresponding to the non-compact semisimple Lie algebras $\gfrak_n$. 
Already this setting is rather involved. We shall point out only 
some complications. One issue is related to the phenomena, that 
unlike the compact forms $U_q(\gfrak_c)$ which are 1-to-1 with $\gfrak_c\,$, 
there are several possible quantum enveloping algebras 
corresponding to a non-compact semisimple Lie algebra $\gfrak_n\,$. 
Moreover, the different procedures to construct $U(\gfrak_n)$ lead to different results, 
cf. e.g., \cite{lyubashenko91,dobrev91,twietmeyer92}. After 
one has selected the deformation of the real form one would 
encounter the next difficulty related to the fact that the 
Casimir operators are not separating the irreducible representations 
well enough. Namely, there are inequivalent irreducible representations 
which share the same values of the Casimir operators. These are 
irreducible representations which are subrepresentations of reducible partially 
equivalent generalized principle series representations. In the 
classical case this phenomenon is explained in, e.g.,  
\cite{dobrevetal77,dobrev88}, and references therein. 

\section*{Acknowledgements}

One of the authours (U.F.) would like to thank M.~Sch\"urmann and M.~Rosso for stimulating discussions related to this work.

\end{document}